\documentclass[11pt, a4paper, twoside]{amsart}
\usepackage{mathrsfs, tabularx}
\usepackage{amsthm, amsmath, amssymb, amscd, float, latexsym, times, epic, eepic}
\usepackage{graphicx,epsfig}

\newtheorem{algo}{Algorithm}[section]

\def \N {{\Bbb{N}}}
\def \R {{\Bbb {R}}}

\renewcommand{\l}{\left}

\newcommand{\norm}[1]{\left\|#1 \right\|}

\def\expec{\mathbb{E}}

\newcommand{\sRV}{S}
\newcommand{\phaseRV}{\Phi}
\newcommand{\noiseRV}{\mathcal{E}}
\newcommand{\Cphase}{C_{\phaseRV}} 
\newcommand{\Cnoise}{C_{\noiseRV}}

\newcommand{\Ltwo}{L^2}			

\def\expec{\mathbb{E}}

\newcommand{\s}{s}					
\renewcommand{\c}{c}				
\renewcommand{\a}{a}

\newcommand{\ngs}{M} 
\newcommand{\nlgs}{{M_\textup{LGS}}} 
\newcommand{\nlay}{L} 

\renewcommand{\l}{\ell}			

\newcommand{\hl}{h_\l} 			

\newcommand{\layer}{\phi}		

\newcommand{\layl}{\layer_{\l}}		
\newcommand{\layfull}{\boldsymbol{\layer}} 

\newcommand{\recfull}{\tilde{\boldsymbol{\layer}}} 

\newcommand{\sfull}{\mathbf{\s}}			

\newcommand{\cfull}{\mathbf{\c}}			

\newcommand{\wav}{\psi}				

\newcommand{\Rec}{\mathbf{R}} 			

\newcommand{\Dfull}{\mathbf{D}}				

\newcommand{\afull}{\mathbf{\a}}	

\newcommand{\outerscale}{K_{out}}
\newcommand{\innerscale}{K_{in}}

\begin{document}

\title[Wavelet methods in multi--conjugate adaptive optics]{Wavelet methods in multi--conjugate adaptive optics}

\author{Tapio Helin and Mykhaylo Yudytskiy}

\maketitle

\begin{abstract}
The next generation ground--based telescopes rely heavily on adaptive optics
for overcoming the limitation of atmospheric turbulence.
In the future adaptive optics modalities, like multi--conjugate adaptive optics (MCAO),
atmospheric tomography is the major mathematical and computational challenge.
In this severely ill-posed problem a fast and stable reconstruction algorithm is needed
that can take into account many real--life phenomena of telescope imaging.
We introduce a novel reconstruction method for the atmospheric tomography problem and
demonstrate its performance and flexibility in the context of MCAO.
Our method is based on using locality properties of compactly supported wavelets,
both in the spatial and frequency domain. The reconstruction in
the atmospheric tomography problem is obtained by solving the Bayesian MAP estimator with a conjugate gradient based algorithm. An accelerated algorithm with preconditioning is also introduced. Numerical performance is demonstrated
on the official end-to-end simulation tool OCTOPUS of European Southern Observatory.
\end{abstract}


\section{Introduction}
\label{sec:intro}


The resolution of an optical imaging system can be limited by several factors. Whereas imperfections of the optical setting can be improved, diffraction always defines a fundamental limit in this regard. 
In ground--based telescope imaging with a small mirror, diffraction is often the dominating effect. However, the influence of the atmospheric
turbulence scales much faster with the increase of the mirror diameter. Already in the next generation of telescopes, called the extremely large telescopes, atmospheric turbulence is the major limiting factor for the angular resolution far beyond the diffraction--limit.
It is a great challenge for science and technology to find ways to achieve diffraction--limited imaging for the future ground--based telescopes.


During the last decades adaptive optics (AO) technology has developed into a powerful remedy for this problem.  Adaptive optics refers to real-time compensation for the distortions in the wavefronts of incoming light due to the atmospheric turbulence. Although the technology involves a number of engineering challenges, the benefits have proven to be fundamental. Adaptive optics correction has been implemented in many major telescope projects, e.g., the Very Large Telescope and the Gran Telescopio Canarias. Moreover, AO is planned as an essential part of all extremely large telescopes.
The work for this paper was largely carried out within a project established towards developing mathematical algorithms for the European Extremely Large Telescope (E-ELT), the future telescope of the European Southern Observatory (ESO).

The mathematical challenges and future prospects for inverse problems field were comprehensively reviewed by Ellerbroek and Vogel in \cite{ElVo09}. With the arrival of next generation implementations of adaptive optics, the severely ill-posed atmospheric tomography problem becomes the crux of mathematical AO research. In this paper we discuss atmospheric tomography in the context of multi--conjugate adaptive optics (MCAO). 

In the classical AO, a single guide star, i.e., a single light source is observed. The aberrations in the incoming light are measured which enables the AO system to correct the cumulative effect of the turbulence towards directions close to the guide star. The correction is performed with a deformable mirror (DM).
Multi--conjugate adaptive optics extends this idea by using several guide stars and multiple deformable mirrors. First, the data obtained from the guide stars is used to solve a reconstruction of the turbulence profile in the atmosphere (see Fig.~\ref{PF_MCAO_fig}). This problem is called atmospheric tomography. Second, having multiple deformable mirrors and the three dimensional reconstruction of the turbulence enable the astronomer to correct for much larger field of view than in the classical implementations. 

The physics of turbulence is an extensive field of study and much is understood about how the turbulence in the atmosphere is formed. Statistical models for turbulence are frequently utilized in the adaptive optics literature by postulating the tomography step as a Bayesian inference problem. In addition, this makes it possible to take into account the statistical nature of the measurement noise. The maximum a posteriori (MAP) estimate is the standard point estimate used to describe the resulting Gaussian posteriori distribution.

In the MCAO related literature, both iterative and non--iterative solution methods have been proposed for the MAP estimator. With the launch of the next generation of extremely large telescopes the dimension of the problem increases rapidly. Even with heavy parallelization, 
the non--iterative methods have a high computational cost and the research in recent years has been inclined towards iterative methods.

The iterative solvers of the MAP estimate are typically based on
conjugated gradient (CG) methods, see \cite{EGV03} and references therein. 
Effort has been put into developing an efficient preconditioner for the problem.
The multigrid preconditioners are investigated in \cite{GEV03,GEV07,GiEl08}, whereas in \cite{VoYa06b, YVE06}, Fourier domain preconditioners have been proposed. 

Especially for iterative methods it is of value to be able to represent the operators in a sparse form. In this regard the Fourier basis is very useful and Fourier-transform based reconstructors have been proposed in \cite{TV01, TLL02, Gavel04}. 
In other typical bases, the forward operator and the inverse covariance of the noise are fast to apply. However, the inverse covariance of the prior is often a full matrix
and thus a sparse approximation is required.
In the approach introduced by Ellerbroek \cite{Ellerbroek02} the turbulence power law is modified in order to achieve a sparse approximation by biharmonics.
Later, a very promising CG based method called the Fractal Iterative Method (FrIM) 
has been developed by Tallon and others in \cite{Tallon_07,  TT10, Tallon_et_al_10, Brunner_et_al_12}. 
There, the inverse covariance is approximated by a factorization that can be applied in ${\mathcal O}(n)$ operations.

We also point out that outside the Bayesian framework the atmospheric tomography problem
has been approached by an algorithm based on the Kaczmarz iteration. The method was introduced by Ramlau and Rosensteiner in \cite{RaRo12}. The authors obtain a very efficient matrix--free 
solver by performing the wavefront reconstruction and the tomography step in two separate problems.
The reconstruction of the atmosphere from incoming wavefronts by the Kaczmarz iteration delivers very promising results in good imaging conditions. The incoming wavefronts are reconstructed from the measurement with an algorithm called the CuReD~\cite{Ros12}.
Several solvers for the wavefront reconstruction exist (see, e.g.,~\cite{Bardsley_et_al_11,PGB02}).

In this paper we suggest a method utilizing compactly supported orthonormal wavelets to represent the atmosphere.
Our method is a CG based iterative method that solves the MAP estimate for the atmospheric tomography problem. We demonstrate a successful performance in low flux imaging conditions, i.e., when only a low number of photons can be measured, and with respect to some practical phenomena that are well--known to limit the reconstruction quality. We have implemented the method with the Daubechies wavelets \cite{Daubechies_88, Cohen_etal_92} in order to have good reconstruction of local details. In addition, the Daubechies wavelets have a useful localization in the frequency domain. This enables us to introduce our key contribution by approximating the inverse covariance with completely diagonal representation. It is shown that such a representation produces an equivalent regularization term
in the MAP estimation problem as indicated by the theory. We point out that this approximation is flexible with respect to choosing a different model for the turbulence power law. In terms of temporal control we rely on an established method called the pseudo--open loop control (POLC) which has been demonstrated to be very robust \cite{Piatrou05}.

In the numerical tests we introduce two variants of the algorithm. In the first setting, we reconstruct more layers of the atmosphere than deformable mirrors using the CG algorithm and optimize the DM shapes accordingly.
In our opinion this demonstrates well the best qualitative performance of the wavelet based method. Moreover, we investigate the stability of the regularization procedure in this setting. The second variant of the method is an accelerated algorithm developed towards achieving the real--time requirements. Here, layers are reconstructed at the altitudes of the deformable mirrors and DM shapes are chosen as the reconstructed layers.
In the accelerated method 
we utilize a modified Jacobi preconditioner, for which we demonstrate fast convergence.   
Numerical simulations are carried out on the OCTOPUS, the official end-to-end simulation tool of ESO. We illustrate the performance of our method in the low flux imaging conditions and compare these results against the matrix--vector multiply (MVM) algorithm, which is the benchmark reconstructor of ESO.

Wavelet methods in adaptive optics have been previously studied in \cite{HaAgBr08, HaAgCoBr10}, however not in the context of MCAO or in the Bayesian scheme. 
In the field of inverse problems wavelets are applied widely (e.g., \cite{SiltanenMueller12, klann_ramlau_reichel}). For an extensive introduction to wavelet basis we refer to \cite{Daub92}.

Notice that atmospheric tomography is a severely ill-posed problem
and is very closely connected to limited angle tomography \cite{Davison}.
Thus, from the general perspective of inverse problems, the theoretical limitations of MCAO are interesting and have been considered in \cite{TV01,TLL02}.
Inverse problems related to waves travelling in random media have been considered, e.g., in the works of Papanicolaou, Bal and Borcea \cite{Borcea_etal, BalPinaud05, Fouque_etal}.

This paper is structured as follows. In Section \ref{sec:math} we discuss the mathematical model for the propagation of light through the atmosphere and how the measurements are obtained. We close the section by explaining the Bayesian paradigm and how the MAP estimate
for the atmospheric tomography is solved. In Section \ref{sec:prior} we introduce the diagonal approximation for the inverse covariance of the turbulence statistics. Section \ref{sec:othereffects} features parts of our method that are essential to MCAO solver, but which have been studied before. Here, we discuss the fitting step and the control algorithm. The concepts of spot elongation and tip--tilt uncertainty are also introduced. These practical phenomena have an essential impact on the noise model. Finally, in Section \ref{sec:numerics} we demonstrate the numerical performance of our method.

\section{Atmospheric tomography}
\label{sec:math}

\subsection{Problem setting}
\label{subsec:problem}

The wind in the atmosphere causes an irregular mixing of warm and cold air. This effect is called the atmospheric turbulence. The fluctuations of the temperature are essentially proportional to the refractive index fluctuations \cite{Ro99} and hence the turbulence affects the propagation of light. With the geometric optics approximation and under appropriate assumptions on the atmosphere, 
the phase of light $\phi$ at the aperture is distorted according to
\begin{equation}
	\label{eq:intg}
	\phi({\bf r}, \boldsymbol\theta) \approx \int_0^H \rho ({\bf r} + \boldsymbol\theta \cdot \xi) d \xi, 
\end{equation}
to a good approximation towards directions $\boldsymbol\theta = (\theta_1, \theta_2, 1)$ close to the zenith. Above, $\rho$ describes the fluctuations of the refractive index, 
${\bf r} = (x,y,0)$ is the location at the aperture, and $H$ is the height of the atmosphere. The approximation \eqref{eq:intg} is derived in \cite{Tatarski1961, Rytov_et_al}.
The challenge in atmospheric tomography is to obtain a good estimate of $\rho$ based on indirect measurements of $\phi(\cdot, \boldsymbol\theta)$ towards  directions $\boldsymbol\theta$ of the guide stars.

The strength of the turbulence at a given altitude varies heavily. However, at the typical telescope sites most of the turbulence is concentrated on certain altitudes. This observation has given rise to a layered atmospheric model, where the refractive index is approximated on a finite number of two-dimensional layers at fixed altitudes. In the simplest example, in ground layer adaptive optics only one layer is considered since the majority of the turbulence strength is located close to the aperture of the telescope (the ground layer). Due to the availability of several deformable mirrors, the implementation of MCAO benefits from a more accurate description of the atmosphere.

\begin{figure}
\centering
\includegraphics[width=0.62\textwidth]{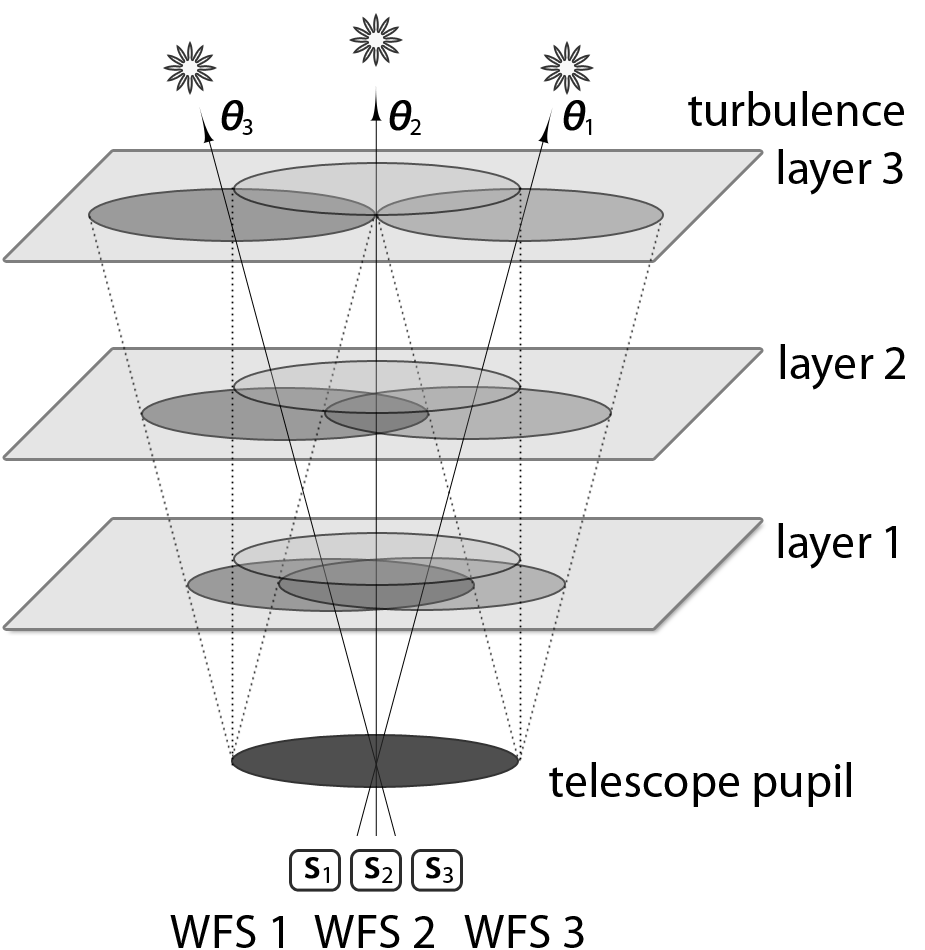}
\caption{In atmospheric tomography the turbulence layers are reconstructed from the wavefront sensor data.}
\label{PF_MCAO_fig}
\end{figure}

Let us consider how the equation \eqref{eq:intg} reduces for a layered atmosphere model.
We denote each modelled layer, located at altitude $\hl$ by $\layl$, $\l = 1, \ldots, \nlay$, 
and by $\boldsymbol\phi = (\phi_1, \ldots, \phi_\nlay)$ a vector representing the atmosphere. 
Assuming geometric propagation, the light arriving from infinity produces incoming wavefronts according to
\begin{equation}
	\label{eq:ngsproj}
	\phi({\bf r}, {\boldsymbol\theta}) = P^\textup{NGS}_{\boldsymbol\theta} \boldsymbol\phi = \sum_{\l = 1}^\nlay \layl({\bf r} + {\boldsymbol\theta} \hl),
\end{equation}
where ${\bf r}$ denotes a point inside the aperture and the vector ${\boldsymbol\theta}$
describes the direction of the guide star (see Fig.~\ref{PF_MCAO_fig}). Here, the projection $P^\textup{NGS}_{\boldsymbol\theta}$
maps the atmosphere to the incoming wavefront from direction ${\boldsymbol\theta}$.
More details on the wave propagation through the layered atmosphere model can be found in~\cite{RoWe96}.

\begin{figure}
\centering
\includegraphics[width=0.62\textwidth]{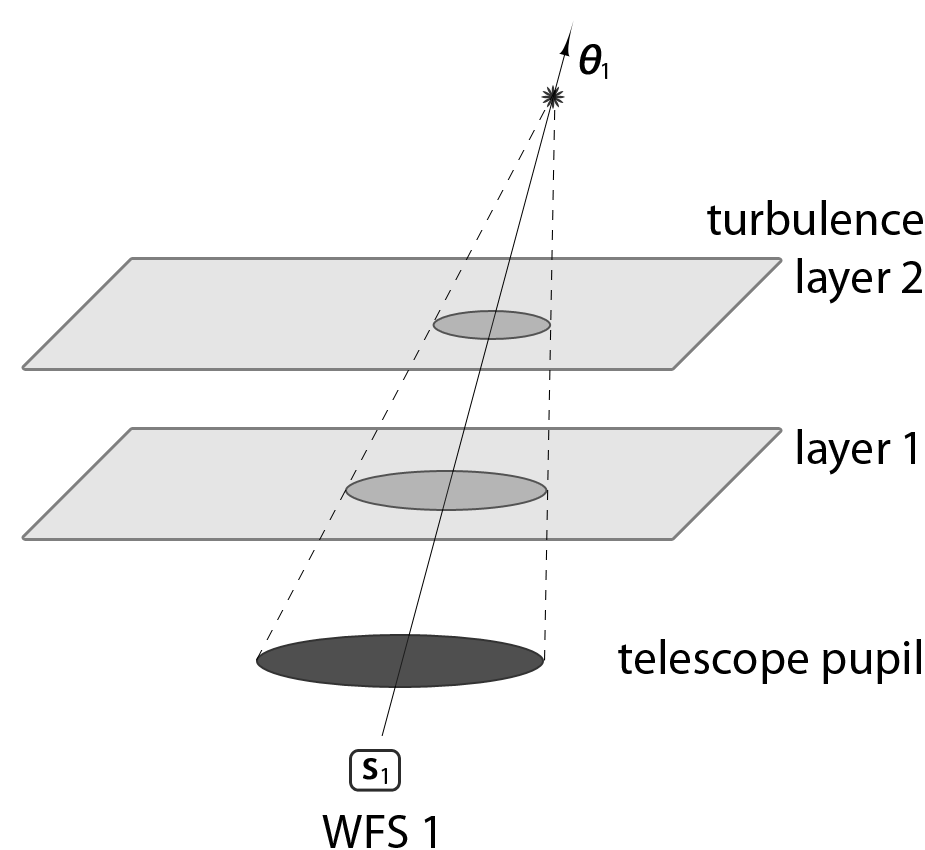}
\caption{
The cone effect.  
Laser guide star light source is fixed above the turbulence layers at some finite altitude. 
Light traveling from the LGS to the telescope pupil passes through smaller areas at higher turbulence layers.
}
\label{PF_FW_CONE_fig}
\end{figure}

However, in practice, there are not enough bright natural guide stars to cover most areas of interest. To overcome this problem, astronomers have developed a technology utilizing lasers, which can generate artificial stars at finite altitude. A laser guide star (LGS) is produced by shooting a powerful laser into the atmosphere where it scatters strongly at certain higher altitudes \cite{Ro99}. Due to the finite altitude, the light arriving to the telescope passes through a cone--like volume in the atmosphere (see Fig.~\ref{PF_FW_CONE_fig}). The corresponding distortion then satisfies
\begin{equation*}
	\phi({\bf r},{\boldsymbol\theta}) = P^\textup{LGS}_{\boldsymbol\theta} \boldsymbol\phi = \sum_{\l = 1}^\nlay \layl\left(\left(1-\frac{\hl}{H}\right){\bf r} + {\boldsymbol\theta} \hl\right),
\end{equation*}
where $H$ denotes the altitude where the laser scatters. Again, $P^\textup{LGS}_{\boldsymbol\theta}$ stands for the projection of the atmosphere to the incoming wavefront with respect to the cone geometry. Whereas an LGS provides a very bright source at directions where no NGS is available, it also introduces some practical limitations. These phenomena, called tip--tilt effect and spot elongation, are described in Section \ref{sec:othereffects}. 

The incoming wavefront can be measured indirectly. A common measurement device
in adaptive optics is the Shack--Hartmann wavefront sensor, which was described in detail in \cite{ElVo09}. Essentially, a Shack--Hartmann sensor measures a quantity proportional to the average gradient of the wavefront on a rectangular grid formed by small lenslets, i.e.,
\begin{equation}
	\label{eq:SH_eq}
	{s}_{ij} = C\int_{\Omega_{ij}} \nabla \phi({\bf r}, \boldsymbol \theta) d{\bf r}
\end{equation}
where ${s}_{ij} = (s^x_{ij}, s^y_{ij}) \in \R^2$ is the measurement and $C$ is a constant. Above, $\Omega_{ij}$ denotes the lenslet also referred to as the subaperture.
Other wavefront sensor modalities exist (e.g., curvature sensor \cite{Ro99}) but are not considered here.

Let us next describe an equation connecting the measurements with the atmosphere for a full MCAO system. The MCAO system we consider will utilize both laser and natural guide stars. Moreover, each guide star here has an individual direction $\boldsymbol \theta$ which we often use to index both the guide star and their corresponding wavefront sensor (WFS). Next, we use $\Gamma$ to denote the measurement operator 
\begin{equation}
	\label{eq:SH_op}
	{\bf s}_{\boldsymbol\theta} = \Gamma_{\boldsymbol\theta} \phi(\cdot, \boldsymbol\theta)
\end{equation}
where ${\bf s}_{\boldsymbol\theta}$ is a vector containing all grid values ${s}_{ij}$ from formula \eqref{eq:SH_eq}. The Shack--Hartmann sensors modelled in equation \eqref{eq:SH_op} can have different resolution and hence $\Gamma_{\boldsymbol \theta}$ is direction-dependent.

In what follows we consider a system that observes $\ngs$ guide stars. Their directions are denoted by $\boldsymbol \theta_m$, where $1\leq m \leq \ngs$. We assume that out of this number, the first $\nlgs$ are laser guide stars. For the rest of the paper we simplify the direction--dependent notations by replacing ${\boldsymbol \theta}_m$ by $m$ whenever no confusion appears.
Now we can write the subproblems for different guide star directions by
\begin{equation}
	\label{eq:atm_tomo_comp}
	{\bf s}_{m} = \Gamma_{m} P^\textup{LGS}_{m}\boldsymbol\phi\quad {\rm and} \quad
	{\bf s}_{m'} = \Gamma_{m'} P^\textup{NGS}_{m'}\boldsymbol\phi
\end{equation}
for $1\leq m \leq \nlgs$ and $\nlgs < m' \leq \ngs$. The full system 
is then described by
\begin{equation}
	\label{eq:atm_tomo}
	{\bf s} = ({\bf s}_{m})_{m=1}^{\ngs} = \mathbf A\boldsymbol\phi,
\end{equation}
where $\mathbf A$ is the concatenation of operators $\Gamma_{m} P^\textup{LGS}_{m}$ and $\Gamma_{m'} P^\textup{NGS}_{m'}$.
Estimating $\boldsymbol\phi$ from a given ${\bf s}$ is called the {atmospheric tomography problem}.

\subsection{Bayesian inference}
\label{subsec:bayes}

The Bayesian inference is a standard approach for solving problem~\eqref{eq:atm_tomo}. This appears natural since information is available about the statistical behavior of the unknown wavefronts and the measurement noise.
The Bayesian paradigm considers problem~\eqref{eq:atm_tomo} in a random setting 
$\sRV = \mathbf A \phaseRV + \noiseRV,$
where $\sRV$, $\phaseRV$ and $\noiseRV$ denote random variables describing the measurement, the incoming turbulence wavefront and the additive noise, respectively.
Given a sample of $\sRV$, i.e., the measurement, the task is to deduce information about the unknown~$\phaseRV$.
Both $\phaseRV$ and $\noiseRV$ are typically modelled as Gaussian random variables. We discuss the distribution of the wavefronts and their physical interpretation below in Section \ref{sec:prior} in more detail. The noise in the Shack--Hartmann measurements is produced by several components \cite{Ro99}.   
However, for the LGS measurements, the effect of spot elongation has a major influence on the noise distribution. We return to the noise covariance and the spot elongation in Section \ref{sec:othereffects}.

Below, we denote the covariance operators of $\phaseRV$ and $\noiseRV$ by $\Cphase$ and $\Cnoise$, respectively. Furthermore, it is assumed that the separate layers are zero centered and uncorrelated. This implies that $\Cphase$ has
a block-diagonal structure
\begin{equation*}
	\Cphase = {\rm diag}\left( C_1, \ldots , C_{\nlay}\right),
\end{equation*}
where $C_{\l}$ denotes the covariance of the layer $\l$.
In the setup given above, the maximum a posteriori estimate can be obtained by solving  
\begin{equation}
\label{REG_REG_map}
\recfull = \operatorname*{argmin}_{\layfull} \left(\|\Cphase^{-1/2}  \layfull \|^2_{2}  + \| \Cnoise^{-1/2} (\sfull -\mathbf A \layfull)\|^2_{2}\right).
\end{equation}
For a general introduction to the Bayesian inverse problems, see \cite{somersalo05}.

\section{Turbulence models and the Bernstein--Jackson equivalence}
\label{sec:prior}


Let us consider for the moment the theory of Gaussian random 
variables in real separable Hilbert spaces. Let $\phi$ be a measurable map from the probability
space $\Omega$ to a Hilbert space $H$. Then $\phi$ is Gaussian if and only if
for all $\rho_1, \ldots, \rho_m \in H$ the mapping $\Omega \ni \omega \mapsto (\langle \phi, \rho_j\rangle)_{j=1}^m$ is a Gaussian random variable in $\R^m$. The distribution of
$\phi$ is determined by the expectation $\expec \phi$ and the covariance operator 
$C_\phi : H \to H$ defined by
\begin{equation*}
	\langle \psi_1, C_\phi \psi_2 \rangle = 
	\expec \left( \langle \phi- \expec \phi, \psi_1 \rangle \langle \phi- \expec \phi, \psi_2 \rangle \right).
\end{equation*}
It is well-known that for any $\psi\in H$ and a linear positive self-adjoint trace class operator $C$ in $H$ with ${\mathcal N}(C) = \{0\}$ there exists a Gaussian random variable $\phi$ in $H$ with mean $\psi$ and covariance $C$ \cite{DPZ}.

Below, we are concerned with zero-centered random variables $\phi$ that have realizations in some Sobolev space $H^s(\R^2)$ with $s>0$ and a covariance operator of the form
\begin{equation}
	\label{eq:def_cphi}
	C_\phi = {\mathcal F}^* M {\mathcal F}.
\end{equation}
Above, ${\mathcal F}$ is the Fourier transform on $\R^2$ 
and $M$ is a multiplication operator $M f(\kappa) = m(\kappa) f(\kappa)$
where $m$ is a positive bounded function. 
With an appropriate decay of $m$ at infinity, the operator $C_\phi$ is trace class and hence $\phi$ is well-defined.
Further, the Schwartz kernel $k_\phi(x,y)$ of the operator $C_\phi$ satisfies
\begin{equation*}
	k_\phi(x,y) = 
	\expec \left( \phi(x)- \expec \phi(x)\right)\left(\phi(y)- \expec \phi(y)\right)
\end{equation*}
in the sense of generalized functions.
We call $k_\phi$ the covariance function of the random process $\phi$.

In the literature of adaptive optics, a turbulence layer $\phi$ is assumed to be isotropic and stationary, i.e., $k_\phi$ depends only on the distance of $x$ and $y$. In particular, $k_\phi$ is
completely characterized by
\begin{equation}
	\label{eq:stationary}
	\tilde k_\phi(z) = k_\phi(x,x+z).
\end{equation}
where $x,z\in\R^2$. Now, it can be shown that if the statistics of $\phi$ satisfy equations \eqref{eq:def_cphi} and \eqref{eq:stationary}, then 
\begin{equation*}
	m(\kappa) = ({\mathcal F} \tilde k_\phi)(\kappa)
\end{equation*}
in the sense of tempered distributions.
The multiplier function $m$ is often referred to as the power spectrum.
The classical model based on the Kolmogorov--Obukhov law of turbulence states that the power spectrum $m$ follows a power law
\begin{equation}
  \label{eq:pwrspec}
  m(\kappa) = C|\kappa|^{-11/3}
\end{equation} 
inside the so-called inertial range $\innerscale \leq |\kappa|\leq \outerscale$ with some constant $C$.
It is not straightforward to expand the power law \eqref{eq:pwrspec} outside the inertial range due to the strong singularity at zero. 
We make a common choice of von Karman model \cite{Ro99, ElVo09} to modify \eqref{eq:pwrspec} by assuming
\begin{equation*}
  m(\kappa) = c_\rho(h) \left(\frac 1{\outerscale^2}+|\kappa|^2\right)^{-11/6},
\end{equation*} 
where $\outerscale$ is the outer scale of the turbulence and $c_\rho(h)$ describes the measure of the optical turbulence strength depending on the altitude. This choice for the power law coincides asymptotically with \eqref{eq:pwrspec} in the high--frequency regime, however, the singularity at zero is removed.
In conclusion, we notice that an equivalence 
\begin{equation}
	\label{eq:cov_approx}
	\|C_{\phi}^{-1/2} f \|_{L^2}^2 
	= \|(\outerscale^{-2}+|\kappa|^2)^{\frac{11}{12}} {\mathcal F} f\|_{L^2}^2
	\simeq \outerscale^{-\frac{11}{3}} \| f\|_{L^2}^2 + \|(-\Delta)^{\frac{11}{12}} f\|_{L^2}^2
\end{equation}
holds in the Cameron--Martin space of $\phi$, i.e., for any $f\in H^{11/6}(\R^2)$.
Here and in what follows, we write $p\simeq q$ if the two pseudo--norms $p$ and $q$ are equivalent.

%

Assume that the wavelets studied here are $r$--regular, i.e., have $r$ vanishing moments and are $r$ times continuously differentiable. For sufficiently large $r$
the last term in \eqref{eq:cov_approx} is equivalent with the expression
\begin{equation}
\label{REG_WAV_norm_equiv}
\| (-\Delta)^{\frac{11}{12}} f \|^2_{L^2} \simeq \sum_{\lambda=1}^\infty 2^{2 \cdot \frac{11}{6} j} | \langle f, \wav_\lambda \rangle |^2,
\end{equation}
where $j$ is the wavelet scale of the wavelet $\wav_\lambda$ with global index $\lambda$. 
The equivalence above is known as the Bernstein--Jackson inequalities \cite{meyer92}.

In the discretized problem, the function $f$ in \eqref{eq:cov_approx} is represented by finite number of wavelet scales. In that case, an equivalent representation for the regularizing term in \eqref{eq:cov_approx} can be produced by a diagonal matrix $D_\ell : \R^{ n_\ell } \to \R^{ n_\ell }$ that satisfies 
\begin{equation*}
	D_\ell = {\rm diag}\left(\frac{1}{c_\rho(h)}\left(\outerscale^{-\frac{11}{3}}+2^{\frac{11}{3}j}\right)\right)_{\lambda=1}^{n_\ell}.
\end{equation*}
Above, $n_{\ell}$ denotes the total number of wavelets for layer $\ell$.
Moreover, we denote $\Dfull = {\rm diag} (D_1, \ldots , D_\nlay)$. Finally, by approximating the
prior covariance of the discretized problem by the ideal model we get
\begin{equation}
\label{REG_WAV_penalty_discr}
\|C_{\Phi}^{-1/2} \layfull \|^2_{(\Ltwo)^\nlay} =
\sum_{\ell=1}^{\nlay} \|C_{\ell}^{-1/2} \phi_{\ell} \|^2_{L^2}
\simeq \sum_{\ell=1}^{\nlay} (D_\ell \mathbf c_\ell, \mathbf c_\ell)_2
= (\Dfull \cfull, \cfull)_{2}
\end{equation}
for $\layfull = (\layl)_{\l = 1}^{\nlay}$ and the wavelet decomposition $\layl =  \sum_{\lambda=1}^{n_\ell} c_{\ell,\lambda} \psi_{\ell,\lambda}$ with respect to the wavelet basis $\{ \psi_{\ell,\lambda} : \lambda = 1,\ldots, n_\ell\}$ of layer $\ell$.
Above, we denote by $\cfull_\ell = (c_{\ell,\lambda})_{\lambda=1}^{n_\ell}$ the vector of wavelet coefficients associated to layer $\ell$, and by $\cfull$ the concatenation of vectors $\cfull_\ell$. 

We point out that in practice the term $\outerscale^{-11/3}$ becomes negligible. Furthermore, the approximation error introduced in \eqref{REG_WAV_penalty_discr} is beyond the scope of this paper. In numerical simulations presented in Section \ref{sec:numerics} we study how different weighting of the regularizing term $(\Dfull \cfull, \cfull)_{2}$ affects the reconstructions obtained by the method.


\section{Other features of MCAO}
\label{sec:othereffects}

\subsection{Fitting step}

In an MCAO system the correction for the wavefront distortions is produced by several deformable mirrors that are conjugated to different altitudes. Hence, a successful mathematical algorithm for MCAO must also deduce optimal mirror shapes for the DMs based on the reconstruction from equation \eqref{eq:atm_tomo}. This subproblem is called the fitting step. Given a sufficient reconstruction of the atmosphere, the fitting step is a well--posed least squares minimization problem and thus classical solution methods provide an efficient reconstruction strategy. We point out that in the ideal fitting one aims
to minimize a functional 
\begin{equation}
	\label{eq:fit_step_ideal}
	{\rm argmin}_{{\bf a}} \expec \left( \int_{{\rm F}} \int_\Omega \left(H_{\boldsymbol\theta} {\bf a} -  P^\textup{NGS}_{\boldsymbol\theta} \boldsymbol \phi\right)^2dx d{\boldsymbol\theta} \right),
\end{equation}
where ${\bf a}$ is the correction profile, $H_{\boldsymbol\theta}$ is the correction towards direction ${\boldsymbol\theta}$ and $P^\textup{NGS}_{\boldsymbol\theta}$ is defined in equation \eqref{eq:ngsproj}. Moreover, $\Omega$ is the aperture domain and $\boldsymbol\theta$ belongs to the field of view $F$.
The problem is typically discretized by choosing a finite set of directions over which the difference in equation \eqref{eq:fit_step_ideal} is averaged. We follow this tradition by
formulating the fitting step as the minimum norm solution to
\begin{equation}
	\label{eq:fit_step_min}
	{\rm argmin}_{{\bf a}}\norm{{\bf H} {\bf a} - {\bf P} \boldsymbol \phi}_2
\end{equation}
where ${\bf H}$ and ${\bf P}$ are concatenations of operators $H_j$ and $P^\textup{NGS}_j$, respectively, towards a finite set of directions ${\boldsymbol\theta}_j$, $j=1, \ldots ,N$ sampled from the field of view.
For more detailed discussion on the fitting step see \cite{ElVo09}.

\subsection{Closed loop control}
\label{subsec:control}

Although in next generation telescopes the DMs are adjusted within milliseconds, the delay between the measurement and the applied DM correction induces an error that needs to be considered. Consequently, a robust temporal control is a fundamental part of the system.

In an MCAO system, the wavefront sensors are located behind the deformable mirrors in the optical path of light.
This is contrary to our assumption on the prior model discussed in the previous section, as the WFS measures the residuals of the incoming wavefronts, instead of the incoming wavefronts themselves.
This mode of operation is called closed loop.


In order to model the physics of turbulence in the prior covariance, we follow here a method called the  pseudo--open loop control (POLC).
The straightforward idea of POLC is to approximate open loop measurements by combing the mirror shapes with the residual data. 
The POLC was introduced to AO in \cite{ElVo03} and further studied in \cite{Gi05_closedstab,Piatrou05}.  It has proven to be stable and robust against large levels of system errors \cite{Piatrou05}.

We rely on a modified POLC, where an integrator is used in the control scheme (see e.g., \cite{Piatrou05}).
We assume that our system has a two time--step delay. Let $t \in  \N, t \ge 2$ denote time--steps. Further, residual Shack--Hartmann data $\sfull$ are measured over the time period $[t-2,t-1)$ for the mirror corrections~$\mathbf a_{t-2}$. We assume that the reconstruction step (computing DM shapes from measurements) consumes one time period, $[t-1,t)$ and the computed mirror shapes $\afull_t$ are applied to the mirrors at time step~$t$.

\begin{algo}{Pseudo--open loop control.}
\label{REG_POLC_algo}
\begin{enumerate}
	\item
		$\sfull^\textup{ol} = \sfull + \widehat{\mathbf A} \afull_{t-2}$
	\item
		$\Delta \afull = \Rec \sfull^\textup{ol} - \afull_{t-2}$
	\item
		$\afull_t = \afull_{t-1} + g \Delta \afull$
\end{enumerate}
\end{algo}

Above, 
$\widehat{\mathbf A}$ maps the mirror shapes $\afull$ to the correction in the measurement space, similar to (\ref{eq:atm_tomo}).
Moreover, the reconstruction operator $\Rec$ maps WFS measurements to DM shapes.
The scalar parameter $0 \le g \le 1$ denotes the gain, which controls the mirror update.

We point out that the POLC is not optimal in the sense of cumulative residual variance. More sophisticated Kalman filter based methods can achieve this \cite{Petit}. The drawback however in such methods is the additional computational load that is a limiting factor, especially for MCAO. Numerically efficient solutions towards this end are an interesting topic for future research.

\subsection{Spot elongation and tip--tilt uncertainty}



Measurements observed from an artificially generated laser guide star suffer from two deficiencies: spot elongation and tip--tilt uncertainty.  These effects were discussed in detail in \cite{ElVo09} and we only briefly state how we correct for these errors in our algorithm.
A successful mathematical model must take these effects into account, as the performance of the AO system would degrade otherwise.

The spot elongation effect occurs due to the physics of scattering at the sodium layer. In practice, a laser guide star is not an ideal point source but rather an extended three dimensional source.
This translates to an elongation of the observed spot on the measurement device, which can be described by a correlation of $x$- and $y$-measurements in each individual subaperture of the Shack--Hartmann WFS.
Our method handles spot elongation following the approach taken in \cite{SpotElong}.

Let us give a brief overview of the noise covariance matrix $C_\mathcal{E}$.
For an MCAO system, which relies on a combination of laser guide star and natural guide star wavefront sensors, the full noise covariance matrix is given as a block-diagonal matrix with respect to the sensors,
\begin{equation*}
C_\mathcal{E} = {\rm diag}\left(\widetilde C_{1}, \ldots, \widetilde C_{\nlgs},
\widetilde C_{\nlgs+1}, \ldots, \widetilde C_{\ngs}\right).
\end{equation*}
Here we associate a noise covariance matrix $\widetilde C_{m}$ for each direction of the guide stars $\boldsymbol \theta_m$, $m = 1, \ldots, \ngs$.
Recall that the first $\nlgs$ directions are associated with sensors observing laser guide stars. The remaining natural guide star sensors are not affected by spot elongation.  
For those directions we assume that the noise is identically and independently distributed in all subapertures with variance $\sigma^2$.  

Also, the spot-elongated measurements are uncorrelated between different subapertures. However, 
for any subaperture the noise in the $x$- and $y$-measurements is correlated and hence the covariance matrix is block diagonal, composed of 2$\times$2 blocks
\begin{equation*}
\widetilde C_{m} = {\rm diag} \left( \widetilde C_{m, 1}, \ldots, \widetilde C_{m, S}\right),
\end{equation*}
where $S$ is the total number of subapertures of the WFS in direction ${\boldsymbol\theta}_m$.
Each block $\widetilde C_{m, i}$, $i = 1, \ldots, S$, can be expressed as
\begin{equation}
\label{eq:elongation_matrix}
\widetilde C_{m, i} = \sigma^2 \left(
I
+ \frac{\tau}{f^2} \boldsymbol\beta \cdot \boldsymbol\beta^\top
\right)
\end{equation}
where $f$ is the full width at half maximum (FWHM) of the non-elongated spots, 
$\boldsymbol\beta$ is the elongation vector and $\sigma^2$ is like above. Moreover, the parameter $0\leq \tau\leq 1$ is used to tune the relative increase of the noise with respect to the elongation.
In consequence, the block-diagonal structure of the covariance matrix implies that applying $C_{{\mathcal E}}$ or its inverse is computationally very cheap. 
For details on deriving $\widetilde C_{m, i}$ we refer to \cite{SpotElong} and the references therein.


The tip--tilt uncertainty is closely related to the uncertainty of the location where the LGS scatters in the atmosphere. Such an error has a large impact on the wavefront sensor measurements. However, it can be shown that the major part of the error is contained in the average $x$- and $y$-derivatives over the whole sensor, i.e., the tip and the tilt of the incoming wavefront.

There are several tip--tilt correction methods that can be applied,
such as a split tomography approach \cite{GiEl08}, a coupled--equation approach \cite{ElVo09}  or a noise--weighted approach \cite{VoYa06b}.
We use a more straightforward approach, in which we remove the incorrect tip--tilt component in the LGS measurements  by modifying equation~(\ref{eq:atm_tomo_comp}) to
\begin{equation}
	\label{eq:removing_tt}
	(I-T) {\bf s}_{m} = (I-T) \Gamma_{m} P^\textup{LGS}_{m}\boldsymbol\phi
\end{equation}
for $m = 1, \ldots, \nlgs$, where $T$ is an orthonormal projection into the tip and tilt components. Other way of stating \eqref{eq:removing_tt} is to say that we use
\begin{equation}
	\label{eq:widehatc}
	\widehat C_{m}^{-1} = (I-T)\widetilde C_{m}^{-1}(I-T)
\end{equation}
as the covariance matrix instead of $\widetilde C_{m}^{-1}$.
Hence this approach neglects more information as, e.g., the noise--weighted approach and relies more on the NGS measurements.
The successful performance of this method is supported by numerical tests.


\section{Numerical implementation}
\label{sec:numerics}

\subsection{Simulation environment and algorithm}

For the simulations, we use the proposed multi-conjugate adaptive optics configuration for the European Extremely Large Telescope.
The telescope gathers light through a circular pupil with diameter of 42~m, of which roughly 28 percent are obstructed.
There are six laser guide stars positioned in a circle with a diameter of 2 arcmins.
To each laser guide star, a Shack--Hartmann wavefront sensor with 84$\times$84 subapertures is assigned.
Moreover, there are three natural guide stars positioned in a circle with a diameter of 8/3 arcmins.
The sensors assigned to the natural guide stars are low resolution Shack--Hartmann sensors (one with 2$\times$2 and two with 1$\times$1 subapertures) for tip--tilt correction.  Further, the E-ELT uses a configuration of three deformable mirrors, located at altitudes 0 km, 4 km and 12.7~km.  The mirrors are modeled by piecewise bilinear functions, with a total number of 9296 degrees of freedom.

We demonstrate the performance of our method on OCTOPUS \cite{OCTOPUS}, the official end-to-end simulation tool of the European Southern Observatory.
The software generates nine frozen layers of the atmosphere located at altitudes between 47 m and 18000 m. The evolution of the atmosphere is simulated by shifting these layers according to their wind directions and speed.

\begin{figure}
\centering
\includegraphics[scale=.85]{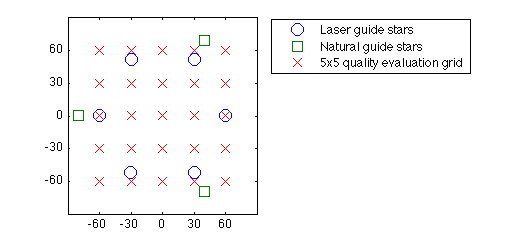}
\caption{Star asterism of the six laser guide stars and three natural guide stars, as well as the 5$\times$5 quality evaluation grid over the field of view (in arcsec).}
\label{grid_lgs_ngs_5x5}
\end{figure}

In the test cases below we simulate one second of evolving atmosphere.  The Shack--Hartmann measurements are read out 500 times per second.  
A two--step delay is observed, as described in Section~\ref{subsec:control}.
The measurements suffer from photon noise and readout noise.
The quality of the reconstruction is evaluated in 25 directions arranged in a 5$\times$5 grid over the field of view, which is a square of 2 arcmins.  As a criteria the long exposure (LE) Strehl ratio \cite{Ro99} is used in K band (for a wavelength of 2200 nm).  The Strehl ratio is a commonly used measure of quality in the astronomical community. Towards directions close to the zenith it can be estimated by the Mar\'echal approximation \cite{Ro99} according to
\begin{equation*}
	s(\boldsymbol\theta) \approx e^{-(2 \pi r({\boldsymbol \theta}) / \lambda)^2},
\end{equation*}
where $s$ is the Strehl ratio, $r({\boldsymbol \theta})$ is the root mean square error in the correction of the incoming wavefront from direction ${\boldsymbol \theta}$, and $\lambda$ is the wavelength.  The long exposure Strehl relates to the average Strehl ratio over the observed timespan.
The star asterism, as well as the 25 evaluation directions are depicted in Fig.~\ref{grid_lgs_ngs_5x5}.

We base our algorithm on the POLC method described in Section~\ref{subsec:control}, where the operator $\mathbf R$ combines the solution operator to the atmospheric tomography problem (\ref{REG_REG_map}) and to the fitting step equation~(\ref{eq:fit_step_min}).  In the special case of reconstructing turbulence layers directly at DM altitudes, we omit the fitting step equation and determine the mirror shapes by the shape of the reconstructed turbulence layers.

The atmospheric tomography problem (\ref{REG_REG_map}) discretized in the wavelet basis is equivalent to solving the linear system of equations
\begin{equation}
	\label{eqn_lse}
	(\widetilde{\mathbf A}^T \Cnoise^{-1} \widetilde{\mathbf A} + \alpha \Dfull) \cfull = \widetilde{\mathbf A}^T \Cnoise^{-1} \sfull,
\end{equation}
where $\widetilde{\mathbf A}$ is the discretization of (\ref{eq:atm_tomo}).
The role of the scalar parameter $\alpha$ is further discussed in Section~\ref{section_test_case3}.
We solve equation (\ref{eqn_lse}) using either the conjugate gradient (CG) method or the preconditioned conjugate gradient (PCG) method with a modified Jacobi preconditioner, discussed below.

In the numerical simulations we utilize the Daubechies wavelet basis. They are a well--known orthogonal wavelet family with compact support \cite{Daubechies_88, Cohen_etal_92}
and have a useful time-frequency localization. For a sufficiently large $n$, the Daubechies~$n$ wavelets are $2$--regular and fulfill the equivalence~\eqref{REG_WAV_norm_equiv}. In order to enhance the spatial localization 
we have chosen to use $n=3$. It is well--known that the Daubechies~3 wavelets
belong to the H\"older space $C^{1+\delta}(\R^2)$ with $\delta\approx 0.0878$~\cite{Daub92}. Even though they are not 2--regular, they form a Riesz basis in $H^2(\R^2)$~\cite{Dahmen95}. 
Based on our numerical tests we believe that this is sufficient in practice.

By formulating the problem in the wavelet domain we gain a significant improvement in terms of convergence of the CG algorithm. This follows from the underlying operators possessing a favorable spectral structure.
Further, the convergence is accelerated by using a modified Jacobi preconditioner, which we discuss in the following.
The operator appearing on the left-hand side of  equation (\ref{eqn_lse}) is given by
\begin{equation*}
	\sum_{m=1}^{\nlgs} (\widetilde A^\textup{LGS}_{m})^T \widehat C_{m}^{-1} \widetilde A^\textup{LGS}_{m}
	+  \sum_{m'=\nlgs+1}^{\ngs} (\widetilde A^\textup{NGS}_{m'})^T \widetilde C_{m'}^{-1} \widetilde A^\textup{NGS}_{m'}
	+ \alpha \Dfull,
\end{equation*}
where $\widehat C_{m}^{-1}$ is defined by equation \eqref{eq:widehatc},
\begin{equation*}
	\widetilde A^\textup{LGS}_{m} = \Gamma_{m} P^\textup{LGS}_{m} W^{-1} \textup{ and } \widetilde A^\textup{NGS}_{m} = \Gamma_{m} P^\textup{NGS}_{m} W^{-1}.
\end{equation*}
Above, $W^{-1}$ is the inverse wavelet transform mapping wavelets to functions and~$\Gamma_m$ is the discretization of the Shack--Hartmann operator according to the Fried geometry (see, e.g., \cite{Ro99}).
Following the discussion of Ellerbroek and Vogel in \cite{ElVo09}, we choose our preconditioner based on only the LGS components, as the low--rank perturbations, corresponding to only a finite number of eigenvalues, do not affect the asymptotic convergence rate of the conjugate gradient algorithm.  Thus, our modified Jacobi preconditioner is
\begin{equation*}
	J = \operatorname*{diag} \left( \sum_{m=1}^{\nlgs} (\widetilde A^\textup{LGS}_{m})^T \widetilde C_{m}^{-1} \widetilde A^\textup{LGS}_{m} \right) + \alpha \Dfull.
\end{equation*}
Finally, we reduce the number of conjugate gradient iterations by choosing the initial guess as the reconstruction in the previous time--step.  This widely used technique for iterative methods in adaptive optics \cite{ElVo09} is known as warm restart.

\subsection{Stability of the regularization} 
\label{section_test_case3}

The diagonal regularization operator in our method was obtained by using the Bernstein--Jackson equivalence. Clearly, the argument applied here does not state explicitly which value for $\alpha$ in \eqref{eqn_lse} is the optimal choice. Also, in this context $\alpha$ can be seen as a regularization parameter. Increasing its value can be considered as stabilization against modeling errors. 
We point out that there are several components of the problem that affect the stability, e.g., the temporal control (gain) and modeling of spot elongation. However, too large value will reduce the quality of the reconstructions. In the following, we demonstrate the performance of the method when $\alpha$ varies.

We study a realistic noise--contaminated situation where the LGSs illuminate 100 photons per subaperture and time--step. Furthermore, the spot elongation and tip--tilt uncertainty for the LGS are simulated.
The NGS tip--tilt sensors observe 500 photons per subaperture and frame.  The readout noise is set to 3 and 5 electrons per pixel for the LGS and NGS sensors, respectively. 

Our algorithm is set as follows. We reconstruct nine layers at the altitudes of the nine simulated turbulence layers using the CG algorithm with 10 iterations.  Further, we solve the fitting step equation (\ref{eq:fit_step_min}) using the CG algorithm with 4 iterations for optimization directions given by the 5$\times$5 evaluation grid.
We choose a gain of $g=0.4$ for the temporal control. The parameter $\tau$ in \eqref{eq:elongation_matrix} was set to $0.8$ for all test cases.

\begin{figure}
	\includegraphics[scale=.85]{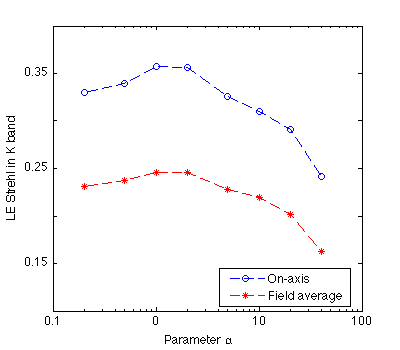}
	\caption{Long exposure Strehl vs.~varying $\alpha$ in equation (\ref{eqn_lse}).}
	\label{figure_test3}
\end{figure}

We run independent simulations with a variable parameter $\alpha$ in (\ref{eqn_lse}) for values $\alpha = 0.2, 0.5, 1, 2, 5, 10, 20$ and $40$.
In Fig.~\ref{figure_test3} we plot the average long exposure Strehl over the field of view (red curve) and the on-axis Strehl (blue curve) for those~$\alpha$. 
The goal of the MCAO system is to obtain the best correction over the field of view, i.e., attain the largest field average Strehl.  However, on-axis Strehl for the zenith-direction is also a quantity of interest.
As can be seen from the plot, the peak on-axis Strehl, as well as the peak field average, is attained with $\alpha = 1$ or~$2$; the difference between the results can be considered negligible.

All higher values of $\alpha$ over-regularize the problem and the performance of the algorithm, although kept stable, decreases.  A low value for $\alpha$ corresponds to an under-regularized problem, and the performance of the method drops.
The importance of the regularizing term $\mathbf D$ in the presence of high photon--noise can clearly be observed.

\subsection{Convergence of the accelerated method}

Here, we demonstrate the convergence properties of our accelerated method.
We run the simulations in the same configuration as in the first test case, where the number of photons per subaperture and time--step for LGS and NGS wavefront sensors are 100 and 500, respectively.  The readout noise is set to 3 electrons per pixel for the LGS sensors and to 5 electrons per pixel for the NGS sensors.

Our accelerated algorithm is set as follows. We reconstruct three layers at the altitudes of the deformable mirrors using the PCG algorithm with the modified Jacobi preconditioner. We utilize $\alpha=1$ in \eqref{eqn_lse} and choose a gain of $g=0.4$ for the temporal control. The parameter $\tau$ in \eqref{eq:elongation_matrix} was set to $0.8$ for all test cases.

\begin{figure}
	\includegraphics[scale=.75]{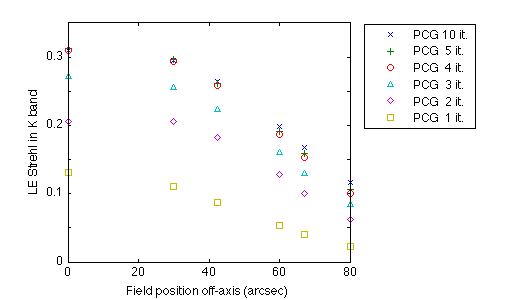}
	\caption{Long exposure Strehl using the PCG method with different number of iterations, averaged over the radial field position.}
	\label{figure_test1}
\end{figure}

In Fig.~\ref{figure_test1} we plot the long exposure Strehl averaged over separation from the zenith after one second of simulated atmospheric propagation.  We run separate simulations for the algorithm with 1, 2, 3, 4, 5 and 10 PCG iterations.  As can be observed, the improvement after iteration 4 is negligible, which indicates that 4~PCG iterations are sufficient for convergence in this configuration.

\subsection{Performance with noisy data}

In this example we consider the performance of our methods with respect to the noise level in the LGS measurements. In other words, we simulate LGSs with different flux between 20 and 200.
We fix the NGS number of photons per subaperture and time--step to 500.
The readout noise is kept as in the previous tests at 3 and 5 electrons per pixel for the LGS and NGS sensors, respectively. 

We compare the performance of our methods with the matrix--vector multiply (MVM) method (see, e.g., \cite{MAD}), which is considered to be the benchmark reconstructor of the ESO.  The MVM is a non--iterative method in which the MAP estimate is discretized using, e.g., the Zernike polynomials. The regularized forward matrix is inverted and applied directly onto the measurements.
The MVM that is presented here reconstructs three layers at DM altitudes, similar to the accelerated wavelet PCG method.

We set the regularization parameter $\alpha = 1$ and the gain $g=0.4$. 
The parameter~$\tau$ in \eqref{eq:elongation_matrix} is tuned for each case separately \cite{SpotElong}.
All tests of the CG method are carried out with 10 iterations for the atmospheric tomography step and 4 iterations for the fitting step; the accelerated method utilizes 4~PCG iterations, which we found to be sufficient above.  

\begin{figure}
	\includegraphics[scale=.85]{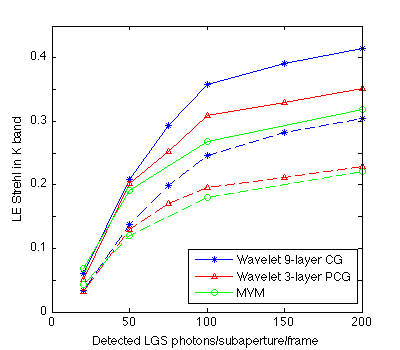}
	\caption{Long exposure Strehl vs.~detected number of photons per subaperture and time--step of LGS sensors.  Solid curve corresponds to the on-axis Strehl; dashed curve to the field average.}
	\label{figure_test2}
\end{figure}

In Fig.~\ref{figure_test2} a comparison in quality of the reconstruction of the three methods is depicted.  Both of the wavelet methods perform better than the MVM in almost all cases; the difference in the 20 photon case can be considered negligible.  We believe the reason for this may be due to a better approximation of the layers using the wavelet basis, as well as the numerical stability of the iterative scheme, as opposed to matrix inversion.

Amongst the two wavelet methods, the approach of reconstructing nine layers followed by a fitting step outperforms the three layer--reconstruction method in quality.  The benefit of the full atmospheric tomography is especially emphasized when more photons are detected by the sensor.  The disadvantage of the nine layer CG method is the higher computational cost over the PCG.


\begin{figure}
	\includegraphics[width=\textwidth]{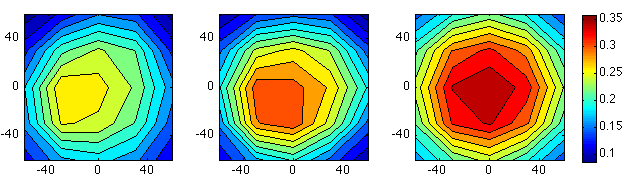}
	\caption{Contour plots of the LE Strehl over the field of view (in arcsec) for MVM (left), 3-layer PCG (middle) and 9-layer CG (right) methods.
	}
	\label{figure_test2_contour}
\end{figure}

To illustrate the difference between the three methods we plot Strehl values in the 25 directions over the field of view for MVM, 3-layer wavelet PCG and 9-layer wavelet CG methods for the $100$ photon case in Fig.~\ref{figure_test2_contour}.

%

\section{Conclusions}

In this paper we have introduced a novel reconstruction method for the atmospheric tomography problem based on wavelets. The theoretical properties of regular wavelets enable us to apply a sparse regularization on the problem that corresponds to utilizing the Kolmogorov turbulence statistics as an a priori model. Here, we have studied the qualitative performance of the method in the context of MCAO.  We derived two variants of the method, concentrating more on quality by solving the full atmospheric reconstruction followed by a fitting step using CG or on speed with the layers--at--DM PCG approach.
We studied the stability of the CG method with respect to the Bernstein--Jackson approximation. Moreover, we demonstrated a fast convergence of the PCG based algorithm. Lastly, we illustrated the quality of the reconstructions in the low--flux regime and showed that the method outperforms the standard reconstruction method, called the MVM, which is used in the ESO simulation platform OCTOPUS. 

We believe that the wavelet method is a very promising algorithm in the field of atmospheric tomography. Fully utilizing the multiscale structure of wavelets can be approached by constructing suitable multigrid preconditioning schemes. Furthermore, the gain in the temporal control can be applied scale--dependently. Together with the careful analysis of the speed of the algorithm we leave these considerations for a future study. We point out that an implementation utilizing the discrete wavelet transform is needed in order to achieve the speed requirements of the E-ELT.

\vspace{1cm}
{\bf Acknowledgements:} This work was done in the framework of the project \emph{Mathematical Algorithms and Software for E-ELT Adaptive Optics}. The project is in cooperation with the European Southern Observatory (ESO) and is funded by the Austrian Federal Ministry of Science and Research. The authors are grateful to the Austrian Adaptive Optics team for support. Moreover, the authors would like to thank Miska Le Louarn, Cl{\'e}mentine B{\'e}chet and Petteri Piiroinen for fruitful discussions. Helin was partly funded by the project ERC-2010 Advanced Grant, 267700.


\clearpage
\addcontentsline{toc}{section}{Bibliography}
\bibliographystyle{plain}
\bibliography{ref_mcao}

\begin{thebibliography}{10}

\bibitem{BalPinaud05}
G.~Bal and O.~Pinaud.
\newblock Time-reversal-based detection in random media.
\newblock {\em Inverse Problems}, 21(5):1593--1619, 2005.

\bibitem{Bardsley_et_al_11}
J.~M. Bardsley, S.~Knepper, and J.~Nagy.
\newblock Structured linear algebra problems in adaptive optics imaging.
\newblock {\em Adv. Comput. Math.}, 35(2-4):103--117, 2011.

\bibitem{SpotElong}
C.~B\'{e}chet, M.~Le Louarn, R.~Clare, M.~Tallon, I.~Tallon-Bosc, and \'{E}.
  Thi\'{e}baut.
\newblock Closed-loop ground layer adaptive optics simulations with elongated
  spots: impact of modeling noise correlations.
\newblock In {\em AO4ELT Proceedings}, 2010.

\bibitem{Borcea_etal}
L.~Borcea, J.~Garnier, G.~Papanicolaou, and C.~Tsogka.
\newblock Enhanced statistical stability in coherent interferometric imaging.
\newblock {\em Inverse Problems}, 27(8):085004, 33, 2011.

\bibitem{Brunner_et_al_12}
E.~Brunner, C.~B\'echet, and M.~Tallon.
\newblock Optimal projection of reconstructed layers onto deformable mirrors
  with fractal iterative method for {AO} tomography.
\newblock In {\em Proc. SPIE, Astronomical Telescopes + Instrumentation}, pages
  84475I--84475I--9, 2012.

\bibitem{Cohen_etal_92}
A.~Cohen, I.~Daubechies, and J.-C. Feauveau.
\newblock Biorthogonal bases of compactly supported wavelets.
\newblock {\em Comm. Pure Appl. Math.}, 45(5):485--560, 1992.

\bibitem{Dahmen95}
W.~Dahmen.
\newblock Multiscale analysis, approximation, and interpolation spaces.
\newblock In {\em Approximation theory {VIII}, {V}ol.\ 2 ({C}ollege {S}tation,
  {TX}, 1995)}, volume~6 of {\em Ser. Approx. Decompos.}, pages 47--88. World
  Sci. Publ., River Edge, NJ, 1995.

\bibitem{Daubechies_88}
I.~Daubechies.
\newblock Orthonormal bases of compactly supported wavelets.
\newblock {\em Comm. Pure Appl. Math.}, 41(7):909--996, 1988.

\bibitem{Daub92}
I.~Daubechies.
\newblock {\em Ten Lectures on Wavelets}, volume~61 of {\em CBMS-NSF Regional
  Conference Series in Applied Mathematics}.
\newblock Soc. for Industrial and Applied Mathematics, Philadelphia, Pa., 1992.

\bibitem{Davison}
M.~Davison.
\newblock The ill-conditioned nature of the limited angle tomography problem.
\newblock {\em SIAM J. Appl. Math.}, 43(2):428, 1983.

\bibitem{Ellerbroek02}
B.~Ellerbroek.
\newblock Efficient computation of minimum-variance wave-front reconstructors
  with sparse matrix techniques.
\newblock {\em J. Opt. Soc. Am.}, 19(9):1803, 2002.

\bibitem{EGV03}
B.~Ellerbroek, L.~Gilles, and C.~R. Vogel.
\newblock Numerical simulations of multiconjugate adaptive optics wavefront
  reconstruction on giant telescopes.
\newblock {\em Appl. Opt.}, 42(24):4811, 2003.

\bibitem{ElVo03}
B.~Ellerbroek and C.~R. Vogel.
\newblock Simulations of closed-loop wavefront reconstruction for
  multiconjugate adaptive optics on giant telescopes.
\newblock {\em Proc. SPIE}, 5169:206--217, 2003.

\bibitem{ElVo09}
B.~Ellerbroek and C.~R. Vogel.
\newblock Inverse problems in astronomical adaptive optics.
\newblock {\em Inverse Problems}, 25:1--37, 2009.

\bibitem{Fouque_etal}
J.-P. Fouque, J.~Garnier, G.~Papanicolaou, and K.~S{\o}lna.
\newblock {\em Wave propagation and time reversal in randomly layered media},
  volume~56 of {\em Stochastic Modelling and Applied Probability}.
\newblock Springer, New York, 2007.

\bibitem{Gavel04}
D.~Gavel.
\newblock Tomography for multiconjugate adaptive optics systems using laser
  guide stars.
\newblock {\em Proc. SPIE}, 5490:1356, 2004.

\bibitem{Gi05_closedstab}
L.~Gilles.
\newblock Closed-loop stability and performance analysis of least squares and
  minimum-variance control algorithms for multiconjugate adaptive optics.
\newblock {\em Appl. Opt.}, 44:993--1002, 2005.

\bibitem{GiEl08}
L.~Gilles and B.~Ellerbroek.
\newblock Split atmospheric tomography using laser and natural guide stars.
\newblock {\em J. Opt. Soc. Am.}, 25(10):2427--35, 2008.

\bibitem{GEV03}
L.~Gilles, B.~Ellerbroek, and C.~R. Vogel.
\newblock Preconditioned conjugate gradient wave-front reconstructors for
  multiconjugate adaptive optics.
\newblock {\em Appl. Opt.}, 42(26):5233, 2003.

\bibitem{GEV07}
L.~Gilles, B.~Ellerbroek, and C.~R. Vogel.
\newblock A comparison of multigrid {V}-cycle versus {F}ourier domain
  preconditioning for laser guide star atmospheric tomography.
\newblock {\em Adaptive Optics: Analysis and Methods/Computational Optical
  Sensing and Imaging/Information Photonics/Signal Recovery and Synthesis
  Topical Meetings}, 2007.

\bibitem{HaAgCoBr10}
P.~J. Hampton, P.~Agathoklis, R.~Conan, and C.~Bradley.
\newblock Closed-loop control of a woofer--tweeter adaptive optics system using
  wavelet-based phase reconstruction.
\newblock {\em J. Opt. Soc. Am. A}, 27(11):A145--A156, Nov. 2010.

\bibitem{HaAgBr08}
P.~J. Hampton and C.~Bradley.
\newblock A new wave-front reconstruction method for adaptive optics systems
  using wavelets.
\newblock {\em Selected Topics in Signal Processing, IEEE Journal of}, 2(5):781
  --792, Oct. 2008.

\bibitem{somersalo05}
J.~Kaipio and E.~Somersalo.
\newblock {\em Statistical and Computational Inverse Problems}, volume 160 of
  {\em Applied Mathematical Sciences}.
\newblock Springer Science+Business Media, Inc, 2005.

\bibitem{klann_ramlau_reichel}
E.~Klann, R.~Ramlau, and L.~Reichel.
\newblock Wavelet-based multilevel methods for linear ill-posed problems.
\newblock {\em BIT Numerical Mathematics}, pages 1--26, 2011.
\newblock 10.1007/s10543-011-0320-x.

\bibitem{OCTOPUS}
M.~Le Louarn, C.~V\'{e}rinaud, V.~Korkiakoski, N.~Hubin, and E.~Marchetti.
\newblock Adaptive optics simulations for the {E}uropean {E}xtremely {L}arge
  {T}elescope.
\newblock In {\em Proc. SPIE 6272, Advances in Adaptive Optics II}, page
  627234, 2006.

\bibitem{MAD}
E.~Marchetti, N.~Hubin, E.~Fedrigo, et~al.
\newblock {MAD} the {ESO} multi-conjugate adaptive optics demonstrator.
\newblock In {\em Proc. SPIE 4839, Adaptive Optical System Technologies II},
  pages 317--328, 2003.

\bibitem{meyer92}
Y.~Meyer and D.~H. Salinger.
\newblock {\em Wavelets and Operators}.
\newblock Cambridge Studies in Advanced Mathematics. Cambridge University
  Press, 1992.

\bibitem{Petit}
C.~Petit, J.-M. Conan, C.~Kulcsár, and H.~F Raynaud.
\newblock Linear quadratic {G}aussian control for adaptive optics and
  multiconjugate adaptive optics: experimental and numerical analysis.
\newblock {\em Journal of the Optical Society of America. A, Optics, image
  science, and vision}, (6):1307–1325, 2009.

\bibitem{Piatrou05}
P.~Piatrou and L.~Gilles.
\newblock Robustness study of the pseudo open-loop controller for
  multiconjugate adaptive optics.
\newblock {\em Appl. Opt.}, 44(6):1003--1010, Feb 2005.

\bibitem{PGB02}
L.~A. Poyneer, D.~T. Gavel, and J.~M. Brase.
\newblock Fast wave-front reconstruction in large adaptive optics systems with
  use of the {F}ourier transform.
\newblock {\em J. Opt. Soc. Am.}, 19(10):2100, 2002.

\bibitem{DPZ}
G.~Da Prato and J.~Zabczyk.
\newblock {\em Stochastic equations in infinite dimensions}, volume~44 of {\em
  Encyclopedia of Mathematics and its Applications}.
\newblock Cambridge University Press, Cambridge, 1992.

\bibitem{RaRo12}
R.~Ramlau and M.~Rosensteiner.
\newblock An efficient solution to the atmospheric turbulence tomography
  problem using {K}aczmarz iteration.
\newblock {\em Inverse Problems}, 28(9):095004, 2012.

\bibitem{Ro99}
F.~Roddier.
\newblock {\em Adaptive Optics in Astronomy}.
\newblock Cambridge University Press, Cambridge, 1999.

\bibitem{RoWe96}
M.~Roggeman and B.~Welsh.
\newblock {\em Imaging Through Turbulence}.
\newblock Cambridge University Press, Cambridge, 1996.

\bibitem{Ros12}
M.~Rosensteiner.
\newblock Wavefront reconstruction for extremely large telescopes via {CuRe}
  with domain decomposition.
\newblock {\em J. Opt. Soc. Am. A}, 29(11):2328--2336, Nov 2012.

\bibitem{Rytov_et_al}
S.~M. Rytov, Y.~A. Kravtsov, and V.~I. Tatarski.
\newblock {\em {Principles of Statistical Radiophysics 4}}.
\newblock Springer--Verlag, Berlin, 1989.

\bibitem{SiltanenMueller12}
S.~Siltanen and J.~L. M\"uller.
\newblock {\em Linear and Nonlinear Inverse Problems with Practical
  Applications}.
\newblock SIAM. 2012.

\bibitem{Tallon_et_al_10}
M.~Tallon, I.~Tallon-Bosc, C.~B\'echet, F.~Momey, M.~Fradin, and \'E.
  Thi\'ebaut.
\newblock Fractal iterative method for fast atmospheric tomography on extremely
  large telescopes.
\newblock In {\em Proc. SPIE 7736, Adaptive Optics Systems II}, pages
  77360X--77360X--10, 2010.

\bibitem{Tallon_07}
M.~Tallon, E.~Thi\'{e}baut, and C.~B\'{e}chet.
\newblock A {F}ractal {I}terative {M}ethod for {F}ast {W}avefront
  {R}econstruction for {E}xtremely {L}arge {T}elescopes.
\newblock In {\em Adaptive Optics: Analysis and Methods/Computational Optical
  Sensing and Imaging/Information Photonics/Signal Recovery and Synthesis
  Topical Meetings on CD-ROM}, page PMA2. Optical Society of America, 2007.

\bibitem{Tatarski1961}
V.~I. Tatarski.
\newblock {\em {Wave propagation in a turbulent medium}}.
\newblock Dover, New York, 1961.

\bibitem{TT10}
E.~Thi{\'e}baut and M.~Tallon.
\newblock Fast minimum variance wavefront reconstruction for extremely large
  telescopes.
\newblock {\em J. Opt. Soc. Am.}, 27(5):1046, 2010.

\bibitem{TLL02}
A.~Tokovinin, M.~Le Louarn, and M.~Sarazin.
\newblock Isoplanatism in a multiconjugate adaptive optics system.
\newblock {\em J. Opt. Soc. Am.}, 17(10):1819, 2002.

\bibitem{TV01}
A.~Tokovinin and E.~Viard.
\newblock Limiting precision of tomographic phase estimation.
\newblock {\em J. Opt. Soc. Am.}, 18(4):873, 2001.

\bibitem{VoYa06b}
C.~R. Vogel and Q.~Yang.
\newblock Fast optimal wavefront reconstruction for multi-conjugate adaptive
  optics using the {F}ourier domain preconditioned conjugate gradient
  algorithm.
\newblock {\em Optics Express}, 14(17), 2006.

\bibitem{YVE06}
Q.~Yang, C.~R. Vogel, and B.~Ellerbroek.
\newblock Fourier domain preconditioned conjugate gradient algorithm for
  atmospheric tomography.
\newblock {\em Appl. Opt.}, 45(21):5281, 2006.

\end{thebibliography}

\end{document}